\documentclass[12pt]{article}
\newcommand{\qed}{\hfill\rule{4pt}{8pt}}

\def\bl{\begin{lm}}\def\el{\end{lm}}
\def\bc{\begin{co}}\def\ec{\end{co}}
\def\bt{\begin{te}}\def\et{\end{te}}
\def\bp{\begin{pr}}\def\ep{\end{pr}}
\def\br{\begin{re}}\def\er{\end{re}}
\def\brs{\begin{res}}\def\ers{\end{res}}
\def\beg{\begin{eg}}\def\eeg{\end{eg}}
\def\bd{\begin{de}}\def\ed{\end{de}}
\def\bdd{\begin{dd}}\def\edd{\end{dd}}
\def\be{\begin{ex}}\def\ee{\end{ex}}
\def\bea{\begin{eqnarray*}}\def\eea{\end{eqnarray*}}

\def\ot{\otimes}

\newtheorem{de}{Definition}[section]
\newtheorem{dd}[de]{}
\newtheorem{lm}[de]{Lemma}
\newtheorem{pr}[de]{Proposition}
\newtheorem{co}[de]{Corollary}
\newtheorem{re}[de]{Remark}
\newtheorem{res}[de]{Remarks}
\newtheorem{ex}[de]{Example}
\newtheorem{te}[de]{Theorem}

\begin{document}
\pagestyle{myheadings}
\title{Lifting of Nichols Algebras of Type $B_2$
\footnote{with an appendix "A generalization of the $q$-binomial
theorem" with Ian Rutherford, Mount Allison University. }}

\author { M. Beattie\thanks{Research partially supported by
NSERC. A visit to U. Bucharest in 1999 was partially supported by
CNCSIS , Grant C12.  She would like to thank the department for
their warm hospitality. } \and S. D\u{a}sc\u{a}lescu
\thanks{Received partial support from CNCSIS, Grants C12 and
199. Thanks to Mount Allison U. for
their hospitality in June 1999. }\\
\and S. Raianu\thanks{On leave from University of Bucharest,
Faculty
of Mathematics. }\\
}
\date{}
\maketitle
\begin{abstract}
We compute liftings of the Nichols algebra of a Yetter-Drinfeld
module of Cartan type $B_2$ subject to the small restriction that
the diagonal elements of the braiding matrix are primitive $n$th
roots of 1 with odd $n\neq 5$. As well, we compute the liftings
of a Nichols algebra of Cartan type $A_2$ if the diagonal elements
of the braiding matrix are cube roots of 1; this case was not
completely  covered in previous work of Andruskiewitsch and
Schneider. We study the problem of when the liftings of a given
Nichols algebra are quasi-isomorphic.  The Appendix (with I.
Rutherford)  contains a generalization of the quantum binomial
formula. This formula was used in the computation of liftings of
type $B_2$ but is also of interest independent of these results.
\end{abstract}

\section{Introduction and preliminaries}
 Let $k$ be an algebraically closed field of characteristic
zero. Several classification results for finite dimensional
pointed Hopf algebras have been obtained in recent years (see
\cite{and} for a survey). The most powerful general method for
classifying such Hopf algebras is the lifting method developed by
N. Andruskiewitsch and H.-J. Schneider. If $A$ is a finite
dimensional pointed Hopf algebra with coradical $k\Gamma$,
$\Gamma$ a group, then there exists a Hopf algebra projection
from $gr(A)$, the associated graded Hopf algebra, to $k\Gamma$,
and this projection splits the inclusion of $k\Gamma$ in $gr(A)$
as the degree 0 component. Then the subalgebra $R$ of
$k\Gamma$-coinvariants of $gr(A)$, called the diagram of $A$, has
 a Hopf algebra structure in the braided category
$^{k\Gamma}_{k\Gamma}{\cal YD}$ of Yetter-Drinfeld modules over
$k\Gamma$. One can also associate to $A$ the Yetter-Drinfeld
module $V$ of primitive elements of $R$, called the infinitesimal
braiding of $A$. The Hopf algebra $gr(A)$ can be reconstructed by
bosonization from $R$, i.e. $gr(A)\simeq R\# k\Gamma$, the
biproduct in the sense of D. Radford or S. Majid. The lifting
procedure consists first in finding all the possible diagrams
$R$, then bosonizing to $gr(A)$, and finally lifting the
information (i.e. presentation by generators and relations) from
$gr(A)$ to $A$.

Assume that $\Gamma$ is a fixed finite abelian group. If $V$ is a
Yetter-Drinfeld module, the Nichols algebra ${\cal B}(V)$ is a
graded Hopf algebra in the category $^{k\Gamma}_{k\Gamma}{\cal
YD}$ with $k1$ as the homogeneous component of degree 0, $V$ as
the homogeneous component of degree 1, and ${\cal B}(V)$ is
generated in degree 1 as an algebra. Nichols algebras were
introduced in \cite{n} (see \cite{ag} for a general presentation
of the construction of and recent developments in Nichols
algebras). Their role in the classification theory for pointed
Hopf algebras was emphasized in \cite{as3}. A fundamental
question is whether the diagram $R$ of $A$ is just the Nichols
algebra ${\cal B}(V)$ of the infinitesimal braiding of $A$. A
positive answer to this question is equivalent to proving the
conjecture that any finite dimensional pointed Hopf algebra is
generated as an algebra by the grouplike elements and the
skew-primitive elements. Up to this conjecture, the lifting
method for classifying finite dimensional pointed Hopf algebras
$A$ with coradical $k\Gamma$ reduces to finding all the
Yetter-Drinfeld modules $V$ such that ${\cal B}(V)$ is finite
dimensional, then  describing  the Nichols algebra ${\cal B}(V)$
by generators and relations for any such $V$, and finally finding
$A$ such that the associated graded Hopf algebra $gr(A)$ is
isomorphic to the biproduct ${\cal B}(V)\# k\Gamma$. Such an $A$
is called a lifting of ${\cal B}(V)\# k\Gamma$.

 A major step in the classification problem was
done in \cite{as3}, where the approach was from the point of view
of Lie theory.  For certain Hopf algebras $A$ (or for any $A$ if
the exponent of $\Gamma$ is prime), the infinitesimal braiding
has a generalized Cartan matrix as an invariant. Then the
dimension of ${\cal B}(V)$ and the structure of this algebra,
reflecting that of $A$, depend on this Cartan matrix and on its
Dynkin diagram. As an example, the lifting method was used in
\cite{as4} to describe liftings of Nichols algebras of Cartan
type $A_2$, and as a consequence classify pointed Hopf algebras
of dimension $p^4$, with $p$ an odd prime . Also the lifting method was
used in \cite{gr} to classify pointed Hopf algebras of dimension
32.

The main aim of this paper is to compute liftings of Nichols
algebras of Cartan type $B_2$. The description of these Nichols
algebras is known (see \cite{as3} and \cite{ringel}). We follow
the general approach that was  used in \cite{as4} for type $A_2$.
The problem of lifting the generators and relations from $gr(A)$
to $A$ has a combinatorial nature, and compared to the  $A_2$ case,
the case of Cartan type $B_2$ requires  more complicated combinatorics.
This is because the structure of the positive roots, which define
a system of generators for the Nichols algebra, is more
complicated in type $B_2$. To deal with these combinatorial
difficulties, we use a generalization of the quantum binomial
formula presented in the Appendix.  In Section 2 we compute the
liftings in type $B_2$. We require that the diagonal elements of
the braiding matrix are primitive $n$-th roots of  odd order not
equal to 5. In fact, in type $A_2$ there was also a case for which
the computation in \cite{as4} failed, more precisely the
case where the diagonal elements of the braiding matrix were
primitive roots of unity of order 3. In Section 3, we show  how
this remaining case can be completed.

 The first
examples of infinite families of nonisomorphic Hopf algebras of
the same dimension were liftings of quantum linear spaces
\cite{as2},\cite{gelaki}, \cite{bdg} or \cite{invent}, and E.
M\"{u}ller's family of nonisomorphic nonpointed Hopf algebras
with nonpointed duals \cite{mueller}.  However, A. Masuoka
\cite{mas} showed that these infinite families consist of Hopf
algebras that are all quasi-isomorphic, i.e. that any element of
the family is a cocycle twist of any other,  or equivalently,
their categories of comodules are monoidally Morita-Takeuchi
equivalent (see \cite{mas} or \cite{peter}). We prove in Section
3 that for $n\neq 5$ and $V$ of type $B_2$ or for $n \neq 3$ and
$V$ of type $A_2$,  any two liftings of ${\cal B}(V)\# k\Gamma$
are quasi-isomorphic.

Recall from any standard text
(such as \cite{k}) the notation
for
$q$-factorials and $q$-binomial coefficients .
Set $(0)_q=0$ and for $n>0,
(n)_q= q^{n-1}+q^{n-2}+\ldots +1$. Set $(0)!_q = 1$ and for $n>0,
(n)!_q = (n)_q(n-1)_q \ldots 1$ .
 Then
 ${n \choose i}_q = \frac{(n)!_q}{(n-i)!_q(i)!_q}$
where $0 \leq i \leq n$. If $i, n $ or $n-i$ is negative,
 then we set ${n \choose i}_q=0$ .

\begin{te}\label{classic} (i). (The q-Pascal identity). For $n \ge k \ge 1,$
\begin{displaymath}\label{qpascal}
{n \choose k}_q
 = {n-1 \choose k-1}_q + q^k{n-1 \choose k}_q
= {n-1 \choose k}_q + q^{n-k}{n-1 \choose k-1}_q.
\end{displaymath}

(ii). (The q-binomial theorem.)
For $x,z $ elements of some $k$-algebra with $zx=qxz$, $q \in k^*$, then
\begin{eqnarray*}\label{qbin}
(x+z)^n = \sum_{i=0}^{n}{n \choose i}_q x^iz^{n-i}.
\end{eqnarray*}  \qed
\end{te}

For $A$ a pointed Hopf algebra with coradical $k\Gamma$, we
denote by $P(A)_{g,h} \mbox{ or } P_{g,h}   $ , if $A$ is clear,
the set $\{x: x \in A, \triangle(x)= g \otimes x + x \otimes h\}.
$ For $\chi \in \hat{\Gamma}, P_{g,h}^\chi = \{x \in P_{g,h}:
lxl^{-1} = \chi(l)x \mbox { for all } l \in \Gamma\}. $

Notation: We write $P_g$ for $P_{g,1}$ and $P_g^\chi$ for
$P^\chi_{g,1}$.

For any coalgebra $C$, $gr(C)$, the graded vector space $C_0
\oplus C_1/C_0 \oplus C_2/C_1 \oplus \ldots$ is a graded
coalgebra.  If $A$ is a pointed Hopf algebra, then $gr(A)$ is a
graded Hopf algebra.

A Yetter-Drinfeld module $V\in{^{k \Gamma}_{k \Gamma}{\cal YD}}$
is a vector space with a left action of $k\Gamma$ and a left
coaction $\delta :V\rightarrow k\Gamma \otimes V$, $\delta
(v)=\sum v_{-1}\otimes v_0$ such that $\delta (hv)=\sum
hv_{-1}h^{-1}\otimes hv_0$ for any $h\in \Gamma$ and $v\in V$.

Throughout, $\Gamma$ will be a fixed finite abelian group and $k$ an
algebraically closed field of characteristic zero.

\section{Liftings of Nichols algebras of type $B_2$}\label{lift}

For $\Gamma$   our fixed finite abelian group ,  let $V$ be a
Yetter-Drinfeld module over $k\Gamma$ of dimension 2. Then $V$
has a basis $\{x_1,x_2\}$ over $k$ such that for $i=1,2$, the
coaction is given by $\delta(x_i)=g_i\otimes x_i$  where $g_i \in
\Gamma$, and the action is given by $g\to x_i = \chi_i(g)x_i$ for
some $\chi_i \in\hat{\Gamma}$,  all $g\in\Gamma$.   The braiding
matrix of $V$ is
\begin{displaymath}
B=\left( \begin{array}{cc}
b_{11} & b_{12}\\
b_{21} & b_{22}
\end{array}\right)
\end{displaymath}
where $b_{ij}=\chi_j(g_i)$. We assume that $V$ has Cartan type
$B_2$, i.e. $b_{ij}b_{ji}=b_{ii}^{a_{ij}}$ for any $i,j$, where
the $a_{ij}$'s are the entries of the Cartan matrix of type $B_2$
\begin{displaymath}
A=\left( \begin{array}{cc}
2 & -1\\
-2 & 2
\end{array}\right).
\end{displaymath}
Thus we have
$$b_{12}b_{21}=b_{11}^{-1},\;\;\; b_{21}b_{12}=b_{22}^{-2},$$
so
\begin{equation}\label{relbij}
b_{12}b_{21}b_{11}=1,\;\; b_{21}b_{12}b_{22}^2=1,\;\;
b_{11}=b_{22}^2
\end{equation}
We assume that $b_{22}=q$ is a primitive root of unity of odd
order $n$ and therefore so is $b_{11}=q^2$. We fix the
Yetter-Drinfeld
module $V$ of type $B_2$ as above throughout all this section.\\

The Nichols algebra ${\cal B}(V)$ has dimension
$n^4$ (see \cite{as3} and \cite{ringel}).
It is presented by generators $x_1, x_2 , z,u$
subject to the relations
\begin{eqnarray*}
z &=& x_2x_1-b_{21}x_1x_2;   \\
u &=& x_2z-b_{21}b_{22}zx_2;   \\
x_1z &=& b_{12}zx_1;    \\
x_2u &=& b_{21}b_{22}^2ux_2;   \\
uz &=& b_{11}b_{21}zu;     \\
x_1 u &=& b_{21}^{-1}b_{12}ux_1+b_{21}^{-1}(b_{22}^{-1}-1)z^2;  \\
x_1^n = x_2^n &=& z^n=u^n=0.
\end{eqnarray*}

Then $H= {\cal B}(V)\# k \Gamma $ has dimension $n^4|\Gamma|$
where $|\Gamma|$ is the order of $\Gamma$. Our goal in this
section is to find all Hopf algebras $A$ such that $gr(A)=H$.

Let $A$ be such a lifting. By \cite[Lemma 5.4]{as2}, we have that
$A$ has coradical $k[\Gamma]$ and $P(A)_h = k(h-1)$ unless
$h=g_i, i=1,2$.  If $h=g_i$ then $P_{g_i}$ has dimension 2 and
$P_{g_i}= k(g_i - 1) \oplus ka_i $ where $ka_i = P_{g_i}^{\chi
_i}$. The image of $a_i$ in $ A_1/A_0$, is just $x_i$. Thus we
have $\triangle(a_i) = g_i \otimes a_i + a_i \otimes 1.$ For $  h
\in \Gamma,\; ha_i = \chi _i(h)a_ih   \mbox { and, in
particular}, g_ja_i = \chi_i(g_j)a_ig_j = b_{ji}a_ig_j.$

Let $\rightharpoonup$ denote the adjoint action and define elements
\begin{equation}\label{defcd}
c = a_2 \rightharpoonup a_1 = a_2a_1 - b_{21}a_1a_2 \;\;\;
\mbox{and} \;\;\; d = a_2 \rightharpoonup c = a_2c -
b_{21}b_{22}ca_2,
\end{equation}

 The following
lemma will be useful in determining the multiplication of the
elements $a_1,a_2, c,d$.

\begin{lm} \label{noskewprim} Let $n$, $\chi_i, g_i$ be as above.
The following assertions hold.
\begin{enumerate}
\item  For $i=1,2$, if $\chi^n_i \neq \epsilon$,
then $P(A)^{\chi_i^n}_{g_i^n} = 0.$
\item  Either $(\chi_1
\chi_2)^n  = \epsilon$ or $P(A)_{(g_1 g_2)^n}^{(\chi_1 \chi_2)^n}
= 0$.
\item  Either $(\chi_1 \chi_2^2)^n  = \epsilon$ or
$P(A)_{(g_1 g_2^2)^n}^{(\chi_1 \chi_2^2)^n} = 0$.
\item  For $n$ different from 5,
$P(A)_{g_1^2g_2}^{\chi_1^2\chi_2} = 0$.
\item   If $n$ is different from 3 and 5, then
$P(A)_{g_1g_2^3}^{\chi_1\chi_2^3} = 0$.
\end{enumerate}
\end{lm}

{\bf Proof.} \begin{enumerate} \item  If $\chi_i^n \neq \epsilon$,
then $P(A)^{\chi_i^n}_{g_i^n} \neq  0$ implies that $\chi_i^n =
\chi_j$ and $g_i^n = g_j$ where $j=1 \mbox{ or } 2$. If $\chi_i^n
= \chi_i$ then $b_{ii}^{n-1}=1$ which is impossible.  If $\chi_i^n
= \chi_j, j \neq i$,
  then $1=b_{ij}$ and $b_{ji}^n = b_{jj}$. But since
$(b_{12}b_{21})^n = 1$, this is impossible.
\item Suppose
$(\chi_1 \chi_2)^n \neq \epsilon$ and $P(A)_{(g_1 g_2)^n}^{(\chi_1
\chi_2)^n} \neq 0$. Then $(\chi_1 \chi_2)^n = \chi_i$ so that
$b_{ii} = b_{ij}^n$. Also $(g_1 g_2)^n = g_i$ so that $b_{ii} =
b_{ji}^n$.  This contradicts $(b_{ij}b_{ji})^n=1.$
\item Suppose
$(\chi_1 \chi_2^2)^n \neq \epsilon$ and $P(A)_{(g_1
g_2^2)^n}^{(\chi_1 \chi_2^2)^n} \neq 0$. If $(\chi_1 \chi_2^2)^n
= \chi_1$ and $(g_1g_2^2)^n =g_1$, then $b^{n-1}_{11} =
b^{-2n}_{12}= b_{21}^{-2n}$ so that $1= b^2_{11}$, which is
impossible. If $(\chi_1 \chi^2_2)^n = \chi_2$ and $(g_1 g^2_2)^n
=g_2$, then $b^{2n-1}_{12} =1= b_{21}^{2n-1}$. But then $1=
(b_{12} b_{21})^{2n-1} = (b_{12}b_{21})^{-1} = b_{11}$, also a
contradiction.
\item For $n \neq 5,  \chi_1^2\chi_2 \neq \epsilon
$. For if $\chi_1^2\chi_2 = \epsilon$, then $b_{11}^2b_{12} =1 =
b_{21}^2b_{22}$. But then $b_{12}= q^{-4} \mbox{  and  } b_{21}^2
= q^{-1}$ so that $ q^{-8}q^{-1} = (b_{12}b_{21})^2 = q^{-4}$ so
that $q^5 =1$, which would imply that $n=5$. Thus if
$P(A)_{g_1^2g_2}^{\chi_1^2\chi_2} \neq 0$, then $ g_1^2g_2 = g_i
\mbox { for } i=1,2.$ But if $g_1^2g_2 = g_1$, then $g_1g_2 = 1$
so that $b_{11}b_{21}=1$ and $b_{12}b_{22}=1$. The first equation
implies $b_{12}=1$ so that the second implies $b_{22}=q=1$, which
is impossible. If $g_1^2g_2 = g_2$ then $g_1^2=1$ so that $q^4=1$
which contradicts the assumption that $n$ is odd.
\item  If
$\chi_1\chi_2^3 = \epsilon$ then $b_{11}b_{12}^3 = 1 =
b_{21}b_{22}^3.$  Then $b_{12}^3 = q^{-2} $ and $b_{21} =
q^{-3}$.  Thus $(b_{12}b_{21})^3 = q^{-11}$ but $(b_{12}b_{21})^3
= q^{-6}$ by (\ref{relbij}). Since $n \neq 5$, this is a
contradiction. Thus if $P(A)_{g_1g_2^3}^{\chi_1\chi_2^3} \neq 0$,
$ P(A)_{g_1g_2^3}^{\chi_1\chi_2^3} = P(A)^{\chi_i}_{g_i} $ for
$i=1,2$.
  If $g_1g_2^3 = g_1$, then $g_2^3 = 1$ and $b_{22}^3 =
q^3 =1$ which contradicts the fact that $n$ is not 3. If
$g_1g_2^3 = g_2$, then $g_1g_2^2 = 1$ so that $b_{11}b_{21}^2 =1
= b_{12}b_{22}^2$. Multiplying the second equation by $b_{21}$
yields $b_{21}=1$ and then the first equation implies $b_{11}=1.$
But then $q=1$ which is impossible.  \qed
\end{enumerate}

Straightforward calculation shows that
\begin{equation}\label{deltac}
\triangle(c) = g_1g_2 \otimes c + c \otimes 1
+ (1-q^{-2})a_2g_1 \otimes a_1.
\end{equation}

Now we compute
\begin{eqnarray*}
\triangle(a_1c)& =& (a_1 \otimes 1 + g_1 \otimes a_1)
(g_1g_2 \otimes c + c \otimes 1 + (1-q^{-2})a_2g_1 \otimes a_1)\\
&=& a_1g_1g_2 \otimes c + a_1c \otimes 1 +(1-q^{-2})a_1a_2g_1 \otimes a_1
 + g_1^2g_2 \otimes a_1c\\
 & +& g_1c \otimes a_1 +
 (1-q^{-2})g_1a_2g_1\otimes a_1^2, \mbox { and }\\
 \triangle(ca_1)& =&
 (g_1g_2 \otimes c + c \otimes 1 + (1-q^{-2})a_2g_1 \otimes a_1)
  (a_1 \otimes 1 + g_1 \otimes a_1)\\
  &=& g_1g_2a_1 \otimes c + ca_1 \otimes 1 +
  (1-q^{-2})a_2g_1a_1 \otimes a_1  \\
 & +& g_1^2g_2 \otimes ca_1
  + cg_1 \otimes a_1 + (1- q^{-2})a_2g_1^2 \otimes a_1^2.
\end{eqnarray*}

Then  using the relations (\ref{relbij}), we see that
  $a_1c - b_{12}ca_1 \in
P(A)_{g_1^2g_2}^{\chi_1^2\chi_2}$ and thus by
Lemma  \ref{noskewprim} (4), if $n \neq 5$,
\begin{equation}\label{a1c}
 a_1c - b_{12}ca_1 = 0.
 \end{equation}

Similarly, using the definition of $d$ and the comultiplication
of $c$ and $a_2$,  we compute
\begin{equation}  \label{deltad}
\triangle (d) = d \otimes 1 + g_1g_2^2 \otimes d +
q(1-q^{-2})a_2g_1g_2 \otimes c +
(1-q^{-2})(1-q^{-1})a_2^2g_1\otimes a_1.
\end{equation}

Further computation shows that $da_2 - b_{12}a_2d$
is $(g_1g_2^3,1)$-primitive and  by Lemma \ref{noskewprim} (5),
if $n \neq 3,5$ then
\begin{equation}\label{da2}
da_2 - b_{12}a_2d=0.
\end{equation}

Now, using (\ref{a1c}) and (\ref{da2}), we compute
\begin{eqnarray*}
da_1 &=& (a_2c - b_{21}b_{22}ca_2)a_1\\
&=& a_2(b_{21}b^2_{22}a_1c) - b_{21}b_{22}c(b_{21}a_1a_2 +c)\\
&=& b_{21}b_{22}^2(b_{21}a_1a_2 + c)c -
b_{21}^2b_{22}(b_{12}^{-1}a_1c)a_2 - b_{21}b_{22}c^2\\
&=& b_{21}^2b_{22}^2a_1(a_2c - b_{21}b_{22}ca_2)
 + b_{21}b_{22}(b_{22} - 1)c^2,
\end{eqnarray*}
so that
\begin{equation}\label{da1}
da_1    = (b_{21}b_{22})^2a_1d + (b_{21}b_{22})(q-1)c^2.
\end{equation}

A similar computation shows that
\begin{equation}\label{cd}
cd = b_{12}dc.
\end{equation}

Now since $a_1$ is $(g_1,1)$-primitive and $q^2=\chi_1(g_1)$ is a
primitive $n$-th root of unity, we have that $a_1^n \in
P(A)_{g_1^n}^{\chi_1^n}$. If $\chi_1^n \neq \epsilon$ then we
must have that  $a_1^n =0$ since $P(A)_{g_1^n}^{\chi_1^n} =0$, by
Lemma \ref{noskewprim} (1), but if $\chi_1^n = \epsilon$, then $a_1^n
= \alpha(g_1^n - 1)$ for some $\alpha \in k$.  By a similar
argument for $a_2^n$ and rescaling  the $a_i$ if necessary , we
have:
\begin{equation}\label{ain}
a_i^n = \mu_i(g_i^n -1), \mbox{ where } \mu_i\in \{0,1 \}, \mu_i=0
\mbox{ if } g_i^n=1 \mbox{ or } \chi_i^n \neq \epsilon.
\end{equation}

\begin{re}\label{bijn}
{\rm If $b_{21}^n \neq 1$, or equivalently $b_{12}^n \neq 1$,
 then $\mu_2 = 0 = \mu_1$.
 For suppose  $b_{21}^n \neq 1$ .  Then
 $\chi_1^n(g_2)=b_{21}^n\neq 1$, so
 $\chi_1^n \neq \epsilon$
 and then $\mu_1 = 0$.  Also since $\chi_2^n(g_1)=b_{12}^n \neq 1$,
 we have that $\chi_2^n \neq \epsilon$, and then $\mu_2 = 0$.}
\end{re}

Now we compute $c^n$ and $d^n$ ; these are the most intricate computations.

By Equation (\ref{deltac}), $\triangle(c) = X + Y + (1-q^{-2}) Z$ where
$X=   g_1g_2 \otimes c, \; Y=   c \otimes 1 ,\; Z=a_2g_1 \otimes a_1$.
Then, since $XY =qYX,  XZ =qZX$, and $q$ is a primitive $n$th root
of unity,  we see that
$\triangle(c^n) = (X+Y+(1-q^{-2}) Z)^n = X^n + (Y+(1-q^{-2}) Z)^n.$
Now
\begin{eqnarray*}
ZY -qYZ &=& a_2g_1c \otimes a_1 - qca_2g_1 \otimes a_1\\
&=& b_{12}b_{11}a_2cg_1 \otimes a_1   - qca_2g_1 \otimes a_1 \\
&=&b_{12}b_{22}^2dg_1 \otimes a_1.
\end{eqnarray*}

Let $T$ denote $dg_1 \otimes a_1$.
Then it is easily checked that $ZT=q^2TZ$ and $TY = q^2YT$ so that
by Remark \ref{serban}  of the Appendix,
$     (Y+(1-q^{-2}) Z)^n  = Y^n +(1-q^{-2})^n  Z^n$.  Thus
\begin{eqnarray*}
\triangle (c^n) & =& (g_1g_2)^n \otimes c^n + c^n \otimes 1
   +  (1-q^{-2})^n(a_2g_1)^n \otimes a_1^n  \\
   &=& (g_1g_2)^n \otimes c^n + c^n \otimes 1
+ (1-q^{-2})^n b_{21}^{-n(n+1)/2}(g_1)^n \mu_2(g_2^n -1)
 \otimes a_1^n\\
   &=& (g_1g_2)^n \otimes c^n + c^n \otimes 1
+ (q^2 -1)^n \mu_2 (g_1)^n (g_2^n -1)
 \otimes a_1^n,\\
 \end{eqnarray*}
since $b_{21}^n=1$ if $\mu_2 \neq 0$ by Remark \ref{bijn}. Let
\begin{equation}\label{upsilon}
\upsilon = c^n + \mu_2(q^2-1)^n a_1^n. \end{equation}Then
$\triangle(\upsilon)$ is
 $$ (g_1g_2)^n \otimes c^n + c^n \otimes 1
+ (q^2 -1)^n \mu_2 (g_1)^n (g_2^n -1)
 \otimes a_1^n + \mu_2(q^2-1)^ng_1^n \otimes a_1^n +
\mu_2(q^2-1)^na_1^n \otimes 1.$$
 It follows that $\upsilon   \in
 P(A)_{(g_1g_2)^n}^{(\chi_1 \chi_2 )^n}$ .  Thus by Lemma
 \ref{noskewprim} (2),
 \begin{equation}\label{cn}
 \upsilon =  c^n + \mu_2(q^2-1)^na_1^n = \lambda (g_1^ng_2^n -1)
\mbox{ where } \lambda=0 \mbox{ if }g_1^ng_2^n
 =1 \mbox{ or }
 (\chi_1 \chi_2 )^n  \neq \epsilon.
\end{equation}

\begin{re} \label{bijnla}
{\rm If $b_{12}^n\neq 1$, then $\lambda =0$. Indeed,
$(\chi_1\chi_2)^n(g_1)=b_{12}^n\neq 1$, which forces
$\lambda =0$.}
\end{re}

We now compute $d^n$.
>From (\ref{deltad}), we see that
 $\triangle(d) = X + Y + bZ +  T$ where
 \begin{eqnarray*}
 X &=&  g_1g_2^2 \otimes d\\
 Y &=&  d \otimes 1\\
 Z &=&  a_2g_1g_2 \otimes c \mbox{ and } b=q(1-q^{-2})\\
 T &=&  (1 - q^{-2})(1-q^{-1}) a_2^2g_1 \otimes a_1\\
 \end{eqnarray*}
 It is easy to check that
 $ XY = q^2 YX, ZY = q^2 YZ, TY = q^2 YT$, so that by
  Theorem (\ref{classic})(ii), the $q$-binomial theorem,
 $$ \triangle (d^n) = (X + bZ +  T)^n + Y^n $$
 and it remains to compute
  $ (X + b Z + T)^n$.
Straightforward computation yields that
 $$XZ = q^2ZX, \mbox { and  } ZT = q^2TZ $$ and we illustrate
 the type of calculation involved by computing $XT$. We have
 \begin{eqnarray*}
XT & =&(1-q^{-1})(1-q^{-2})(g_1g_2^2 \otimes d)(a_2^2g_1 \otimes a_1)\\
 &=& (1-q^{-1})(1-q^{-2})b_{12}^2b_{22}^4a_2^2g_1^2g_2^2 \otimes da_1\\
 &=& (1-q^{-1})(1-q^{-2})b_{12}^2b_{22}^4a_2^2g_1^2g_2^2
  \otimes (b_{21}^2b_{22}^2a_1d + b_{21}b_{22}(q-1)c^2)\\
 &=& (1-q^{-1})(1-q^{-2})( q^2a_2^2g_1^2g_2^2 \otimes a_1d
 +    (q-1)b_{12}q^3a_2^2g_1^2g_2^2  \otimes c^2).
 \end{eqnarray*}
 Since $Z^2 = b_{12}qa_2^2g_1^2g_2^2 \otimes c^2$, we have
 \begin{equation}\label{computext}
 XT = q^2TX + (1-q^{-2})(1-q^{-1})(q-1)q^2Z^2.
 \end{equation}

 Using Theorem \ref{main} in the Appendix, we see that
$$ (X + bZ + T)^n = X^n + \nu(n)Z^n + T^n   $$
 and by Corollary \ref{conu},
  $\nu(n) = \alpha^n + \beta^n$ where $\alpha = q-1 , \beta= 1-q^{-1}$
 are the solutions of the equation $Y^2 - q(1-q^{-2})Y +
 (1-q^{-2})(1-q^{-1})(q-1)q^2/(q^2 -1) = 0$.  Thus
 $\nu(n) = (q-1)^n + (1-q^{-1})^n  = 2(q-1)^n$ and we have
 \begin{eqnarray*}
 \triangle(d^n)&=& (g_1g_2^2)^n \otimes d^n + d^n \otimes 1
 + 2(q-1)^n(a_2g_1g_2)^n \otimes c^n +
  (q^2-1)^n(q-1)^n(a_2^2g_1)^n \otimes a_1^n \\
 & = & (g_1g_2^2)^n \otimes d^n + d^n \otimes 1
  +  2(q-1)^n  a_2^ng_1^ng_2^n \otimes c^n
 + (q^2 -1)^n(q-1)^n a_2^{2n}g_1^n \otimes a_1^n,
  \end{eqnarray*}
where $(a_2g_1g_2)^n = (b_{22} b_{12})^{n(n-1)/2}a_2^ng_1^ng_2^n$
and $b_{12}^{n(n-1)/2}= 1$ if $a_2^n \neq 0$ by Remark \ref{bijn}.
Similarly $(a_2^2g_1)^n = a_2^{2n}g_1^n$. Let
\begin{equation}\label{omega}
\omega = d^n + 2(q-1)^n \mu_2 c^n +(q^2 -1)^n
(q-1)^n \mu_2^2 a_1^n \end{equation}
 and then
\begin{eqnarray*}
\triangle (\omega)& =& (g_1 g_2^2)^n \otimes d^n + d^n \otimes 1 + 2(q-1)^n \mu_2
(g_2^n -1) (g_1g_2)^n \otimes c^n\\
&+& (q^2-1)^n (q-1)^n \mu_2^2 (g_2^n -1)^2 g_1^n \otimes a_1^n \\
&+& 2(q-1)^n \mu_2 g^n_1 g_2^n \otimes c^n + 2(q-1)^n \mu_2 c^n
\otimes 1 + 2(q-1)^n \mu_2^2 (q^2-1)^n (g_2^n -1) g_1^n \otimes a_1^n\\
&+& (q^2 -1)^n (q-1)^n \mu_2^2 a_1^n \otimes 1 + (q^2-1)^n (q-1)^n
\mu_2^2 g^n_1 \otimes a_1^n\\
& =& (g_1 g_2^2)^n \otimes [ d^n + 2(q-1)^n \mu_2 c^n
+ (q^2-1)^n (q-1)^n \mu_2^2 a^n_1]\\
&+& [d^n +2(q-1)^n \mu_2 c^n +(q^2 -1)^n (q-1)^n \mu^2_2 a_1^n ] \otimes 1 \\
&-& 2 (q-1)^n \mu_2 g_1^n g_2^n \otimes c^n
+ (q^2 -1)^n (q-1)^n \mu_2^2 (-2g_2^n +1)g_1^n \otimes a_1^n\\
&+& 2(q-1)^n \mu_2 (g_1g_2)^n \otimes c^n
+   2(q-1)^n (q^2 -1)^n \mu^2_2  (g^n_1) (g_2^n -1) \otimes a_1^n\\
&+& (q^2 -1)^n (q-1)^n \mu_2^2 g_1^n \otimes a_1^n \\
&=& (g_1g_2^2)^n \otimes \omega + \omega \otimes 1.
\end{eqnarray*}

Thus $\omega$ is $((g_1g_2^2)^n,1)-$ primitive and so, by
Lemma \ref{noskewprim} (3), we have
\begin{eqnarray}\label{dn}
 \omega = d^n + 2(q-1)^n \mu_2 c^n +(q^2 -1)^n (q-1)^n
\mu_2^2 a_1^n  &=& \gamma ((g_1g_2^2)^n -1) \\
 \mbox{ where } \gamma  = 0
 \mbox { if } (\chi_1\chi_2^2)^n \neq \epsilon &\mbox{ or }&
(g_1g_2^2)^n =1.    \nonumber
\end{eqnarray}

\begin{re} \label{bijnga}
{\rm If $b_{21}^n\neq 1$, then $\gamma =0$. Indeed,
$(\chi_1\chi_2^2)^n(g_2)=b_{21}^n\neq 1$, so
$(\chi_1\chi_2^2)^n\neq \epsilon$, forcing $\gamma =0$.}
\end{re}

We find a characterization for the liftings of Nichols algebras of
type $B_2$   similar to that in \cite{as4} for type $A_2$.
The following lemma from \cite{as4} will be useful.

\begin{lm} \label{aslemma} \cite[Lemma 3.4 (i)]{as4}
Let $X,Y,Z$ be elements in a $k$-algebra, $\alpha, \beta$ scalars
in $k$ and $n$ a natural number. If $YX = \alpha XY +Z$
and $ZX = \beta XZ$, then
$$
YX^n = \alpha^n X^n Y +(\sum^{n-1}_{i=0} (\alpha^i \beta^{n-1-i}) X^{n-1} Z$$
and if $\alpha \neq \beta$ and $\alpha^n = \beta^n$, then for $n>1$,
$$YX^n = \alpha^n X^n Y.
$$
 \qed
\end{lm}

Define a Hopf algebra $U^+$ in the category $^{k\Gamma}_{k\Gamma}
{\cal YD}$ by $U^+: = k < x_1, x_2 , z,u | {\cal N} >$, where
${\cal N}$ is the set of the first six relations defined at the
beginning of this section, namely:
\begin{eqnarray}
z&=& x_2x_1-b_{21}x_1x_2 \label{defz} \\
u&=& x_2z-b_{21}b_{22}zx_2 \label{defu} \\
x_1z&=&b_{12}zx_1  \label{x1z} \\
x_2u&=&b_{21}b_{22}^2ux_2  \label{x2u}\\
uz &= &b_{11}b_{21}zu   \label{uz} \\
x_1u&=&b_{21}^{-1}b_{12}ux_1+b_{21}^{-1}(b_{22}^{-1}-1)z^2  \label{x1u} ,
\end{eqnarray}
and the comultiplication, action and coaction are defined such that
$x_1, x_2$ are primitive and $h \cdot x_i = \chi_i (h) x_i,
\delta (x_i) = g_i \otimes x_i$.
To see that $U^+$ is well defined, we note that if we make the free
algebra $F$ generated by $x_1$ and $x_2$ into a braided Hopf algebra
with $x_1$ and $x_2$ primitives, then $\Delta _F({\cal N})\subseteq
{\cal N}\ot F+F\ot {\cal N}$, and this induces a braided Hopf algebra
structure on $U^+$.
$U^+$ has a PBW basis $\{\; x_2^iu^jz^rx_1^s\; |\; i,j,r,s\geq 0\; \}$.
This follows from the fact that $U^+$ can be constructed from
$k[x_1]$ by adjoining $z,u,x_2$ by iterated Ore extensions
defined by the relations (\ref{defz})-(\ref{x1u}).

We define $U$  to be the Radford biproduct $U^+ \# k \Gamma$.

\begin{te}\label{quotient}
Let $\mu _1,\mu _2\in \{ 0,1\}$ and $\lambda , \gamma \in k$ such that
\begin{eqnarray}
\mu_i &=& 0
\mbox { if } g_i^n=1 \mbox { or
} \chi_i^n \neq \epsilon; \label{mui} \\
 \lambda &=& 0 \mbox{ if } g_1^ng_2^n=1 \mbox{ or }
\chi_1^n\chi_2^n \neq \epsilon;\\
 \gamma &=& 0 \mbox { if } (\chi_1\chi_2^2)^n \neq
\epsilon \mbox{ or } (g_1g_2^2)^n =1.
\end{eqnarray}
Then the two-sided ideal $J$ of $U$ generated by the elements
$$y_i := x_i^n - \mu_i(g_i^n -1), \;\; i=1,2; $$
$$ v:= z^n + \mu_2 (q-1)^n x_1^n - \lambda(g_1^ng_2^n -1);
 $$
$$ w:= u^n + 2(q-1)^n \mu_2 z^n +(q^2 -1)^n (q-1)^n \mu_2^2 x_1^n  -
\gamma ((g_1g_2^2)^n -1) ,$$
is a Hopf ideal of $U$. Moreover
$A=A(\Gamma ,V,(\mu _i)_i,\lambda ,\gamma )=U/J$
is a pointed Hopf algebra of dimension $n^4|\Gamma |$ with
coradical $k\Gamma$, and $gr(A)\simeq {\cal B}(V)\# k\Gamma$, where
$V$ is our fixed Yetter-Drinfeld module of type $B_2$.
\end{te}
{\bf Proof.}  Since, by  the arguments preceding (\ref{ain}),
(\ref{cn}) and (\ref{dn}), $J$ is generated by skew-primitive
elements, $J$ is a Hopf ideal. Next we verify some commutation
relations. We have
\begin{eqnarray*}
x_2 x^n_1 &=& b^n_{21} x^n_1 x_2 \mbox{ by Lemma \ref{aslemma}
with }Y =x_2, X=x_1, Z=z, \alpha =b_{21}, \beta= b_{21} b^2_{22};\\
zx_1^n &=& b_{21}^n x^n_1 z \mbox{ by } (\ref{x1z}); \\
ux_1^n &=& (x_2 z - b_{21} b_{22} zx_2)x_1^n = (b_{21}^2)^n x^n_1 u;\\
zx_2^n &=& b^n_{12} x^n_2z \mbox{ by Lemma \ref{aslemma}
with $Y =z, X=x_2, Z =u, \alpha= b_{12} b_{22}, \beta= b_{12}$};\\
ux_2^n &=& b^n_{12} x^n_2 u \mbox{ by } (\ref{x2u}); \\
x_1z^n &=& b^n_{12} z^n x_1 \mbox{ by }(\ref{x1z}); \\
x_2z^n &=& b^n_{21} z^n x_2 \mbox{ by Lemma \ref{aslemma}
with $Y =x_2, X=z, Z=u, \alpha = b_{21} b_{22}, \beta= b_{11} b_{21}$};\\
uz^n &=& b^n_{21} z^n u  \mbox{ by }(\ref{uz});\\
x_1 u^n &=& b_{12}^{2n} u^n x_1
\mbox{ by Lemma \ref{aslemma} with } Y=x_1, X=u, Z= z^2,
\alpha = b^2_{12} b^2_{22}, \beta = b^2_{12}; \\
x_2 u^n &=& b^n_{21} u^n x_2 \mbox{ by } (\ref{x2u}); \\
zu^n &=& b^n_{12} u^n z \mbox{ by } (\ref{uz}).
\end{eqnarray*}

It remains to find the commutation between $x_1$ and  $ x_2^n.$

Now  $x_1 x_2^n= b_{12} b^2_{22} (x_2x_1 -z) x_2^{n-1}$, and
$zx_2^t =(b_{12} b_{22})^t x_2^t z- b^t_{12} ( \sum^{t-1}_{i=0}
b^i_{22})b_{22} x_2^{t-1} u$ by Lemma \ref{aslemma}.

Thus $x_1x_2^n = b_{12}b^2_{22} x_2 x_1 x_2^{n-1} -(
b_{12}b^2_{22})(b_{12}b_{22})^{n-1} x_2^{n-1} z +(b_{12}b_{22}^2)
b_{12}^{n-1} ( \sum^{n-2}_{i=0} b^i_{22}) b_{22} x_2^{n-2} u.$

Then $$x_1x_2^n =(b_{12} b^2_{22})^n x_2^n x_1 + \alpha x_2^{n-1}
z + \beta x_2^{n-2}u$$ and we show that $\alpha = \beta =0$.

It is easy to see that
$$\alpha = - b_{12}^n \sum^{n}_{i=1} b^{2i}_{22} b_{22}^{n-i} =
- b^n_{12} \sum_{i=1}^{n} b^i_{22} =0,$$ and that
\begin{eqnarray*}
\beta &=& b^n_{12} b^3_{22}\left \{ \sum^{n-2}_{i=0} b^i_{22}
+ b^2_{22} \sum^{n-3}_{i=0} b^i_{22} + b^4_{22} \sum^{n-4} _{i=0} b^i_{22}
+ \ldots + b_{22}^{2(n-2)} (1)\right \}\\
&=&b^n_{12} b^3_{22}\left \{ (q^{n-1}-1) +q^2(q^{n-2}-1)+\ldots
+q^{2(n-2)}(q-1)\right \} / (q-1). \end{eqnarray*}
% &=&b^n_{12}
% b^3_{22}\left \{ (q^{n-2}-1)/(q-1)+q^{-1}-(q^{2(n-1)}-1)/
% (q^2-1)\right \} \\
% &=&b^n_{12} b^3_{22}((q^{-2}-1)(q+1)+q^{-1}(q^2-1)-q^{-2}+1)/(q^2-1)\\
% &=&0
% \end{eqnarray*}
Now the expression in brackets is just
\begin{eqnarray*}
&& (q^{n-1} + q^2q^{n-2} + q^4q^{n-3} + \ldots + q^{2n-3}) - ( 1 +
q^2 + q^4 + \ldots + q^{2(n-2)})\\
&=& q^{-1}(q^{n-1} - 1)/(q-1) - (q^{2(n-1)} - 1)/(q^2 -1),
\end{eqnarray*}
and putting these expressions over a common denominator, we see
that this is $0$ and so $\beta = 0$.

We have proved that
$$x_1x_2^n=b_{12}^nx_2^nx_1$$
We show now that $J$ is the right ideal generated by $y_i , \;
i=1,2;\;
  v;
\;  w$.\\
Assume first that $b_{12}^n=1$. Then for $h\in \Gamma$ and
$i=1,2$, we have $$hy_i = \chi _i^n(h)y_i h+ (\chi _i^n(h)-1)\mu
_i(g_i^n-1)h$$ and by (\ref{mui}), we always have $(\chi
_i^n(h)-1)\mu _i=0$. Also
$$x_2y_1=y_1x_2 \mbox{  and  }
 x_1 y_1=y_1x_1.$$ Similar
computation for the other generators shows that $J$ is the right
ideal as well as the two-sided ideal generated by $y_1,y_2,v,w$. \\
If $b_{12}^n\neq 1$, then by Remarks \ref{bijn}, \ref{bijnla} and
\ref{bijnga}, we must have $\mu_1=\mu_2=\lambda =\gamma =0$ and
then $J$ is the two-sided ideal generated by
$x_1^n,x_2^n,z^n,u^n$. The commutation relations show immediately
that $J$ is the right ideal generated by these elements.

We prove now that no non-zero linear combination of the elements
$gx_2^iu^jz^rx_1^s$, with $g\in \Gamma$, $0\leq i,j,r,s \leq
n-1$, lies in $J$. This will imply that the dimension of $A=U/J$
is $n^4$ and also that $J\cap k\Gamma =0$, so $k\Gamma$ embeds in
$A$. To show this we proceed as in the proof of \cite[Proposition
1.10]{bdg}. Assume that
$$
\sum _{g,i,j,r,s}\alpha _{g,i,j,r,s}gx_2^iu^jz^rx_1^s = \sum
_{i=1,2}y_i f_i +
 vf_3
+ wf_4
$$
for some $f_1,f_2,f_3,f_4\in U$ and some scalars $\alpha _{g,i,j,r,s}$,
not all equal to zero. The commutation relations show that
$U$ is a free module with basis $\{ \; x_2^iu^jz^rx_1^s\; |
0\leq i,j,r,s\leq n-1\; \}$ over the subalgebra
$B$ of $U$ generated by $\Gamma$ and $x_1^n,x_2^n,z^n,u^n$. If we
write $f_i, 1\leq i\leq 4$ in terms of this basis, we see that there
exist some $F_1,F_2,F_3,F_4\in B$ such that
 $$
\sum _{i=1,2} y_iF_i +
  vF_3
+ wF_4 \in k\Gamma -\{ 0\}
 $$
Clearly $B$ is isomorphic as an algebra to an Ore extension $R$
obtained from $k\Gamma$ by adjoining the indeterminates
$Y_2,Y_3,Y_4,Y_1$ (identified with $x_2^n$, $u^n$, $z^n$, $x_1^n$)
in that order via Ore extensions with zero derivations. This
shows that the relation
\begin{eqnarray*}
\sum _{i=1,2}(Y_i - \mu_i(g_i^n -1))q_i +
 (Y_4 + \mu_2 (q-1)^n Y_1 - \lambda(g_1^ng_2^n -1))q_3 &&\\
+(Y_3 + 2(q-1)^n \mu_2 Y_4 +(q^2 -1)^n (q-1)^n \mu_2^2 Y_1  -
\gamma ((g_1g_2^2)^n -1))q_4&\in & k\Gamma -\{ 0\}
\end{eqnarray*}
holds in $R$ for some $q_i, 1\leq i\leq 4$. The universal property
for Ore extensions (see \cite[Lemma 1.1]{bdg}) shows that there
exists an algebra morphism $\theta :R\rightarrow k\Gamma$ acting
as identity on $\Gamma$ and such that $\theta (Y_i)=\mu
_i(g_i^n-1)$ for $1\leq i\leq 2$, $\theta (Y_4)=-\mu _1\mu
_2(q-1)^n(g_1^n-1)+\lambda (g_1^ng_2^n-1)$, and $\theta
(Y_3)=-2(q-1)^n\mu _2(-\mu _1\mu _2(q-1)^n(g_1^n-1)+ \lambda
(g_1^ng_2^n-1))-(q^2-1)^n(q-1)^n\mu _1\mu _2^2(g_1^n-1) +\gamma
((g_1g_2^2)^n-1)$. Then applying $\theta$ to the above equation
we obtain that $0\in k\Gamma -\{ 0\}$, a contradiction.

We have thus proved that the dimension of $A$ is $n^4|\Gamma |$ and
that $k\Gamma$ embeds in $A$. Since $A$ is generated
by $\Gamma$ and the skew-primitive elements $x_1,x_2$,
we see that $A$ is pointed and the coradical of $A$ is $k\Gamma$.

For the last claim, we consider the algebra morphism $\phi :U
\rightarrow gr(A)$ which takes $x_i$ to the image of $x_i$ modulo $k\Gamma$
in the homogeneous component of degree 1 of $gr(A)$. Since
$A$ is generated as an algebra by $\Gamma$ and $x_1,x_2$, this algebra
morphism is surjective. On the other hand, since $x_i^n\in k\Gamma$
in $A$, we have that $x_i^n=0$ in $gr(A)$. Similarly, regarding the
images of $z$, respectively $u$, in $gr(A)$, they have degrees 2,
respectively 3, and we also get that $z^n=0$ and $u^n=0$ in $gr(A)$.
Therefore $\phi$ induces a surjective algebra morphism $\psi$
from $U/(x_1^n,x_2^n,z^n,u^n)\simeq {\cal B}(V)\# k\Gamma$ to
$gr(A)$, and this morphism must be an isomorphism because of the
dimensions. Obviously $\psi$ is also a coalgebra morphism, and this ends the
proof.
\qed

\begin{te}   \label{b2lift}
Let $A$ be a pointed Hopf algebra with coradical $k\Gamma$ and such
that $gr(A)\simeq {\cal B}(V)\# k\Gamma$, where $V$ is our fixed
Yetter-Drinfeld module of type $B_2$, such that $n$ is odd and
$n\neq 5$. Then
$A\simeq A(\Gamma ,V,(\mu _i)_i,\lambda ,\gamma )$
for some $g_i,\chi _i,\mu _i,\lambda,\gamma$ as in Theorem \ref{quotient}.
\end{te}
{\bf Proof.}
Suppose first that $n\neq 3$.
We have shown that relations
(\ref{a1c}), (\ref{da2}), (\ref{da1}), (\ref{cd}) hold in $A$,
so there exists a Hopf algebra morphism $\phi :
A(\Gamma ,V,(\mu _i)_i,\lambda ,\gamma )
\rightarrow A$ which takes $x_i$ to $a_i$ for $i=1,2$. By
\cite[Lemma 2.2]{as2} we have that $A$ is generated as an algebra
by $\Gamma$, $a_1$ and $a_2$, so $\phi$ is surjective. The
dimension implies now that $\phi$ is an isomorphism.

Now suppose that $n=3$; in this case relation (\ref{da2}) has not
been verified.  Since $n \neq 5$, by the proof of Lemma
\ref{noskewprim} (5), $\chi_1\chi_2^3 \neq \epsilon$ so that
 $da_2 - b_{12} a_2d
\in P^{\chi_1 \chi^3_2}_{g_1 g_2^3}$ means $g_1 g_2^3 = g_i$ and
$\chi_1 \chi_2^3 = \chi_i$ for $i$ = 1 or 2. If $i=2$, then the
argument is the same as in Lemma \ref{noskewprim}. If $i=1$, then
$g_2^3=1$ and $\chi_2^3 = \epsilon$. But then $a_2^3=0$ by
(\ref{ain}) and thus $a_2^3 \rightharpoonup a_1 = 0 = a_2
\rightharpoonup d = a_2d - b_{21}b_{22}^2da_2.$ Relation
(\ref{da2}) follows from (\ref{relbij}). \qed

\section{Quasi-isomorphism of liftings }

Recall that   Hopf algebras $A$ and $B$ are quasi-isomorphic if
one is a cocycle twist of the other  and this implies that their
categories of comodules are monoidally Morita-Takeuchi equivalent
(see \cite{mas} or \cite{peter}). If one of $A$ or $B$ is pointed
or finite dimensional,  then the converse holds. If $A$ and $B$
are quasi-isomorphic, we write $A \sim B$.

As well as the infinite families of nonisomorphic Hopf algebras
of the same dimension   obtained by lifting quantum linear
spaces that were mentioned in the introduction, such infinite
families
  can also be
easily constructed from liftings of ${\cal B}(V) \# k\Gamma$,
where $V \in {^{k\Gamma}_{k\Gamma} {\cal YD}}$ is of type $A_2$ or
$B_2$. Recall that $V$
 of type $A_2$ means that $V = kx_1 \oplus kx_2$
and there exist $g_1,g_2\in\Gamma$ and
$\chi_1,\chi_2\in\hat{\Gamma}$ such that for all $g\in\Gamma$ and
$i=1,2$, we have
$$g\to x_i = \chi_i(g)x_i,\; \mbox {and}\;\;\delta(x_i)=g_i\otimes x_i.$$
Also for $b_{ij} = \chi_j(g_i)$ as in Section \ref{lift},
\begin{equation}\label{relbija2} b_{12}b_{21}b_{11}=1 =
b_{21}b_{12}b_{22},\;\; b_{11}=b_{22}=q, \mbox{ where $q$ is a
primitive $n$th root of unity, }
\end{equation}
so that the Cartan matrix $A$ determined from the braiding matrix
$B=(b_{ij})$   is
\begin{displaymath}
A=\left( \begin{array}{cc}
2 & -1\\
-1 & 2
\end{array}\right) .
\end{displaymath}
The liftings of ${\cal B}(V) \# k\Gamma$ where $V \in
{^{k\Gamma}_{k\Gamma} {\cal YD}}$ is of type $A_2$ were determined in
\cite{as4} where $n>3$ or $\Gamma$ is cyclic of order 3.

For $x_1,x_2$ as above, a Hopf algebra $U^+ \in {^{k\Gamma}_{k\Gamma}
{\cal YD}}$, is defined in \cite[Definition 3.5]{as4}  by
$$U^+ = k<x_1,z,x_2| z= x_1x_2 - b_{12}x_2x_1, \; zx_1 = b_{21}x_1z, \;
x_2z = b_{21}zx_2>,$$
where the $x_i$ are primitive elements.
 Then $U$ is defined to be the Radford
biproduct $U^+ \# k\Gamma$.

\begin{te} \label{lifta2}\cite[Theorems 3.6 and 3.7.]{as4}
\ Let $A$ be a lifting of ${\cal B}(V) \# k\Gamma$ for $V$ of
type $A_2$ as described above.  If $n$ is an odd integer
greater than 3 or if $\Gamma$ is cyclic of order 3, then $A \cong
U/J$, where $J$ is the Hopf ideal of $U$ generated by the
skew-primitives
$$x_i^n - \mu_i(g_i^n -1) \mbox{ where }
\mu_i \in \{0,1\} \mbox{ and } \mu_i = 0 \mbox { if } g_i^n=1
\mbox { or } \chi_i^n \neq \epsilon,$$
$$ z^n + \mu_1 (q-1)^nx_2^n - \lambda(g_1^ng_2^n -1)
\mbox{ where } \lambda=0 \mbox{ if } g_1^ng_2^n=1 \mbox{ or }
\chi_1^n\chi_2^n \neq \epsilon. $$
\end{te}

\begin{ex}\label{infinitea2}
{\rm (cf \cite[Section 3]{as4}) Let $\Gamma = <g>$ be cyclic of order
49. Let $q$ be a primitive 7-th root of unity. Let $\chi \in
\hat{\Gamma}$ be defined by $\chi (g) =q$. Define $g_1=g, g_2 =
g^4, \chi_1 = \chi, \chi_2 = \chi^2$. Let $V= kx_1 \oplus kx_2$
where $x_i \in V_{g_i}^{\chi_i}$ . Then the matrix
$$B = \left [ \begin{array}{cc} b_{11} & b_{12} \\ b_{21} & b_{22}
\end{array} \right ] = \left [ \begin{array}{cc} q & q^2 \\ q^4 & q^8 =q
\end{array} \right ]$$
and the Cartan matrix $A$ is $\left [ \begin{array}{rr} 2 & -1 \\
-1 &2
\end{array} \right ]$ so $V$ is of type $A_2.$

Let $A(\lambda)$ be the Hopf algebra $U / < x_1^7 -( 1-g^7),
x_2^7 - (1-(g^4)^7), z^7 -(q-1)^7(1-(g_4)^7) - \lambda
(1-(gg^4)^7) >$ as in Theorem \ref{lifta2}.

For $\lambda \neq \omega$, $A(\lambda) \not \cong A(\omega)$. For
suppose $f: A(\lambda) \to A(\omega)$ is a Hopf algebra
isomorphism. Then $f(g) =g$ since if $f(g) = g^4$, then $f(g^4) =
g^{16} \neq g$. Thus $f(x_1) = \alpha y_1 + \delta (1-g^7)$,
where $y_i$ is the counterpart of $x_i$ in $A(\omega).$
Commutation with $g$ shows that $\delta =0$ and  $\alpha^7 =1$.
Similarly $f(x_2) = \beta y_2$ where $\beta^7 =1$. Then $f(z) =
\alpha \beta (y_1 y_2 - b_{12} y_2 y_1) = \alpha \beta \omega$,
and $(\alpha \beta)^7 \omega^7 = \omega^7 = (q-1)^7 (1- g^{28}) -
\lambda (1- g^{35}) = (q-1)^7 (1- g^{28}) - \omega (1-g^{35}),$
so $\lambda = \omega$. } \qed
\end{ex}

We now describe the liftings for the remaining case $n=3$. Let
${\cal U}^+$ be the free algebra in the  indeterminates $x_1$ and
$x_2$.  This is a Hopf algebra in the category
${^{k\Gamma}_{k\Gamma} {\cal YD}}$ by taking the $x_i$'s to be
primitive elements where $\delta(x_i) = g_i \otimes x_i$ and $h
\rightharpoonup x_i = \chi_i (h) x_i$ for all $h \in \Gamma$.
Denote $z= x_1x_2 - qx_2x_1$. We define ${\cal U}$ as the Radford
biproduct ${\cal U}^+ \# k\Gamma$.

\begin{pr}\label{a2nequal3}Suppose $V \in {^{k\Gamma}_{k\Gamma} {\cal YD}}$
is of type $A_2$ and $n=3$. Then any lifting of ${\cal B}(V) \#
k\Gamma$ is isomorphic to ${\cal U}/J$ where ${\cal U}$ is as
defined above and $J$ is the Hopf ideal of ${\cal U}$ generated
by:
$$x_i^3 - \mu_i(g_i^3 -1) \mbox{ where }
\mu_i \in \{0,1\} \mbox{ and } \mu_i = 0 \mbox { if } g_i^3=1
\mbox { or } \chi_i^3 \neq \epsilon;$$
$$zx_i - q^ix_iz - \gamma_i(g_i^2g_j -1), i \in \{1,2\}, i \neq j,
\gamma_i = 0 \mbox{ if } g_i^2g_j = 1 \mbox{ or } \chi_i^2\chi_j
\neq \epsilon; $$
$$ z^3 + \mu_1 (q-1)^3x_2^3 + (1-q)\gamma_1(zx_2 - q^2 x_2z)
- \lambda(g_1^3g_2^3 -1)
\mbox{ where } \lambda=0 \mbox{ if } g_1^3g_2^3=1 \mbox{ or }
\chi_1^3\chi_2^3 \neq \epsilon. $$
\end{pr}

{\bf Proof.} Let $A$ be a lifting of ${\cal B}(V) \#
k\Gamma$.
Let $a_i \in A$ be the lifting of $x_i \in V$ as in the
$B_2$ case and as in    \cite{as4},
$a_i^3 = \mu_i(g_i^3 -1) \mbox{ where } \mu_i \in \{0,1\} \mbox{
and } \mu_i = 0 \mbox { if } g_i^3=1 \mbox { or } \chi_i^3 \neq
\epsilon$. For $c= a_1a_2 - b_{12}a_2a_1$ as in \cite{as4},  then
it is shown in \cite[Lemma 3.1]{as4} that $ca_1- b_{21} a_1c \in
P^{\chi_1^2 \chi_2}_{g_1^2 g_2}$ and $a_2c- b_{21} ca_2 \in
P^{\chi_1 \chi_2^2}_{g_1 g_2^2}$. If $\chi_i^2\chi_j \neq
\epsilon$ for $i,j \in \{1,2 \}$, or if $\Gamma$ is cyclic of order 3,
then we are in the situation of
Theorem \ref{lifta2}. But for $n=3$ and $\Gamma$ not cyclic of order 3,
we could have
$\chi_i^2\chi_j = \epsilon $; then the matrix $B = (b_{ij})
=\left (
\begin{array}{cc} q & q\\q&q
\end{array} \right )$ with $q^3 =1$. From now on, we assume this
is the case. Then
\begin{equation}\label{cai}
ca_i - q^i a_ic = \gamma_i(g_i^2 g_j-1)   \mbox{ for some
}\gamma_i \in k \mbox{ with }\gamma _i= 0 \mbox{ if }g_i^2 g_j=1
\mbox{ or } \chi_i^2\chi_j \neq \epsilon.
\end{equation}
Now as in the calculation of (\ref{deltac}), we have
$$\triangle(c) = g_1g_2 \otimes c +
c \otimes 1 + (1-q^2)a_1g_2 \otimes a_2 = X + Y + (1-q^2)Z .$$
Then we have the following commutation relations:
\begin{eqnarray*}
XY &=& qYX; \\
XZ &=& qZX + \gamma_2q^2T \mbox{ where } T=a_1g_1g_2^2 \otimes
(g_1g_2^2 -1);\\
YZ &=& q^2ZY + \gamma_1S \mbox{ where } S= (g_1^2g_2 -1)g_2
\otimes a_2;\\
XT &=& q^2TX ;\\
ZT &=& qTZ;\\
YT &=& qTY + \gamma_1 W \mbox{ where } W = (g_1^2g_2-1)g_1g_2^2
\otimes (g_1g_2^2 - 1);\\
XS &=& q^2SX + \gamma_2W;\\
YS &=& qSY;\\
ZS &=& q^2SZ;\\
TS &=& q^2ST.
\end{eqnarray*}

Direct  computation shows that
\begin{eqnarray*}
 & & \triangle(c^3) = (X + Y + (1-q^2)Z)^3\\
&=& (g_1g_2)^3 \otimes c^3 + c^3 \otimes 1 + (q-1)^3a_1^3g_2^3
\otimes a_2^3 + (1-q)\gamma_1\gamma_2g_1g_2^2(g_1^2g_2 -1)
\otimes (g_1g_2^2 -1).
\end{eqnarray*}

Let  $\upsilon = c^3 + \mu_1 (q-1)^3 a_2^3 + (1-q)\gamma_1
\gamma_2(g_1g_2^2 -1)$. Then
\begin{eqnarray*}
\triangle (\upsilon) &=& (g_1g_2)^3 \otimes c^3 + c^3 \otimes 1
+(q-1)^3 \mu_1 (g_1^3 -1) g^3_2 \otimes a^3_2\\
&+& (1-q) \gamma_1 \gamma_2 g_1 g_2^2 (g_1^2 g_2-1) \otimes (g_1
g_2^2 -1) + \mu_1 (q-1)^3 a_2^3 \otimes 1 + \mu_1 (q-1)^3 g_2^3
\otimes a_2^3\\
&+&(1-q) \gamma_1 \gamma_2 (g_1 g_2^2 -1) \otimes 1 +(1-q)
\gamma_1 \gamma_2 g_1 g_2^2 \otimes (g_1 g^2_2 -1)\\
&=& (g_1 g_2)^3 \otimes [c^3 +(q-1)^3 \mu_1 a_2^3 +(1-q) \gamma_1
\gamma_2 (g_1 g^2_2 -1)]\\
&+& [c^3 + \mu_1 (q-1)^3 a_2^3 +(1-q) \gamma_1 \gamma_2 (g_1 g_2^2
-1)] \otimes 1\\
&-&  \mu_1 (q-1)^3 g_2^3 \otimes a_2^3 -(1-q) \gamma_1 \gamma_2
g_1g^2_2 \otimes (g_1 g^2_2 -1) + \mu_1 (q-1)^3 g_2^3 \otimes
a^3_2\\
&+& (1-q) \gamma_1 \gamma_2 g_1
g^2_2 \otimes (g_1 g^2_2 -1)\\
&=& (g_1 g_2)^3 \otimes \upsilon +\upsilon \otimes 1,
\end{eqnarray*}
and thus $\upsilon  \in   P^{\chi_1^3 \chi_2^3}_{g_1^3 g_2^3} .$
If $\chi_1^3 \chi_2^3 \neq \epsilon$, then $ \chi_1^3 \chi_2^3  =
\chi_i$ for $i=1,2$, yielding $q=1$, which is a contradiction.
Therefore, $\upsilon = \lambda(g_1^3g_2^3 - 1)$ for some $\lambda
\in k$. Now an argument
similar to the one in Theorem \ref{quotient} shows that
the elements $ha_1^ic^ja_2^k, h \in \Gamma, 0 \leq i,j,l,k \leq
2$ are a basis for $A$, and the same argument as in Theorem
\ref{b2lift} completes the proof.
 \qed

Let $A(\Gamma,V,\mu_1,\mu_2,\lambda,\gamma_1,\gamma_2)$ denote
the Hopf algebra ${\cal U}/J$ where ${\cal U}$ is the Hopf
algebra defined just before Proposition \ref{a2nequal3} and $J$
is the Hopf ideal generated by the skew primitives:
$$x_i^n - \mu_i(g_i^n -1) \mbox{ where } i\in \{1,2\},
\mu_i \in \{0,1\} \mbox{ and } \mu_i = 0 \mbox { if } g_i^n=1
\mbox { or } \chi_i^n \neq \epsilon;$$
$$zx_1 - b_{21}x_1z - \gamma_1(g_1^2g_2 -1), \gamma_1 = 0
\mbox{ if } g_1^2g_2 = 1 \mbox{ or } \chi_1^2\chi_2 \neq
\epsilon; $$
$$zx_2 - b_{21}^{-1}x_2z - \gamma_2(g_2^2g_1 -1),
\gamma_2 = 0 \mbox{ if } g_2^2g_1 = 1 \mbox{ or } \chi_2^2\chi_1
\neq \epsilon; $$
$$ z^n + \mu_1 (q-1)^nx_2^n + (1-q)\gamma_1(zx_2 - b_{21}^{-1} x_2z)
- \lambda(g_1^ng_2^n -1) \mbox{ where } \lambda=0 \mbox{ if }
g_1^ng_2^n=1 \mbox{ or } \chi_1^n\chi_2^n \neq \epsilon. $$ Then
by Theorem \ref{lifta2} and Proposition \ref{a2nequal3}, all
liftings of Nichols algebras of type $A_2$ are of this form. Now
we show that all
$A(\Gamma,V,\mu_1,\mu_2,\lambda,\gamma_1,\gamma_2)$ with the same
$\gamma_i$ are quasi-isomorphic. We write $H \sim H'$ if the Hopf
algebras $H$ and $H'$ are quasi-isomorphic. We use the key theorem
from \cite{mas} together with comments from \cite{maspersonal}.

Recall that for $K$ a Hopf algebra, the set $Alg(K,k)$ is a group
under the convolution product with the inverse to $\psi \in
Alg(K,k)$ being given by $\psi \circ S$ where $S$ is the antipode
of $K$.  The left action of $Alg(K,k)$ on $K$ is given by $\psi x
= (Id_K \otimes \psi )\triangle(x)$ and the right action by $x
\psi = (\psi \otimes Id_K)\triangle (x)$.  Two Hopf ideals $I,J$
in a Hopf algebra $K$ are said to be conjugate if there is an
algebra map $\psi$ from $K$ to $k$ such that $J=\psi I \psi^{-1}$
.. Also if $K$ is a subHopf algebra of a Hopf algebra $H$ and $J$
is a Hopf ideal of $K$, then $(J)$ will denote the Hopf ideal in
$H$ generated by $J$.

\begin{te}\cite[Theorem 2]{mas} \cite{maspersonal} \label{masuoka}
Suppose that $K$ is a Hopf subalgebra of a Hopf algebra $H$. Let
$I,J$ be Hopf ideals of $K$.  If there is an algebra map $\psi$
from $K$ to $k$ such that $J=\psi I \psi^{-1}$  and $H/ (\psi I)$
is nonzero, then $H/ (\psi I)$ is an $(H/(I), H/(J))$-biGalois
object and so    the quotient Hopf algebras $H/(I),H/(J)$ by the
Hopf ideals $(I),(J)$ in $H$ generated by $I,J$ are monoidally
Morita-Takeuchi equivalent.
\end{te}

In the application of Masuoka's theorem, the following lemma will be useful.

\begin{lm} \label{conjugate} Let $K$ be a Hopf algebra
containing  $(g_i, 1)$-primitives $x_i, i =1, \ldots, t$. Let $J$
be the Hopf ideal of $K$ generated by the $x_i$ and let $L$ be the
Hopf ideal generated by $x_i - \lambda_i (g_i -1), i=1, \ldots,
t$. Let $\psi$ be an algebra map from $K$ to $k$ such that $\psi
(x_i) = \lambda_i$ and $\psi (h) =1$ for $h$ grouplike. Then $J$
and $L$ are conjugate ideals in $K$.
\end{lm}

{\bf Proof.} Since $S(x_i) = - g_i^{-1} x_i, \psi (S(x_i)) = -
\lambda_i$. Thus
$$\psi^{-1} x_i \psi = ( \psi \otimes Id \otimes \psi^{-1})
(\triangle^2 x_i)
= \psi (g_i) g_i (- \lambda_i) + \psi (g_i) x_i + \psi (x_i)
= x_i - \lambda_i (g_i -1),$$
and $\psi^{-1} J \psi =L$. \qed

\begin{te} \label{quasia2} For  $V$ of type $A_2$, and
$A=A(\Gamma,V,\mu_1,\mu_2,\lambda,\gamma_1,\gamma_2)$ a lifting
of ${\cal B} (V) \# k \Gamma$, then $A$ is quasi-isomorphic to
any other lifting  $A(\Gamma,V,\mu_1',\mu_2',\lambda',
\gamma_1,\gamma_2)$.  If $\gamma_1=0$, then $A \sim
A(\Gamma,V,\mu_1',\mu_2',\lambda', 0,\gamma_2')$ and if $\gamma_2
=0$, then  $A \sim A(\Gamma,V,\mu_1',\mu_2',\lambda',
 \gamma_1',0).$ In particular, if $n>3$ or $\Gamma$ is cyclic of
order 3, then all liftings of ${\cal B} (V) \# k \Gamma$ are
quasi-isomorphic.
\end{te}

{\bf Proof.} First note that if $b_{ij}^n \neq 1$, for $i \neq j$,
then for $n > 3$ by an argument similar to Remark \ref{bijn},
$b_{ji}^n \neq 1$ and  the only possible lifting is the trivial
one, $A(\Gamma,V,0,0,0,0,0)$. For $n=3$ this follows from the
proof of Proposition \ref{a2nequal3}.

Now assume $b_{21}^n = b_{12}^n =1$. Since $n$ is odd, then also
$b_{21}^{n(n-1)/2} =1$. First, we show that for a given
$\gamma_1,\gamma_2,\mu_1$,
$A(\Gamma,V,\mu_1, 0,0,\gamma_1,\gamma_2)$ and $A=A(\Gamma,V,\mu_1,
\mu_2,\lambda,\gamma_1,\gamma_2)$ are quasi-isomorphic for any
$\mu_2 \in \{ 0,1\}$ and any $\lambda$.

Let $M_{\mu_1} =  {\cal U}/< x^n_1 - \mu_1 ( g_1^n-1),\; zx_1 -
b_{21}x_1z - \gamma_1(g_1^2g_2 -1),\; zx_2 - b_{12}b_{22}x_2z -
\gamma_2(g_1g_2^2-1)>$ where ${\cal U}$ is the Hopf algebra
defined just before  Proposition \ref{a2nequal3}. Note that
$\upsilon = z^n +\mu_1(q-1)^n  x^n_2 +
(1-q)\gamma_1\gamma_2(g_1g_2^2 -1)$ is $(g_1^n
g_2^n,1)$-primitive in $M_{\mu_1}$ . If $\gamma_1 = \gamma_2 =0$,
this follows from the proof of Theorem 3.1 \cite[Theorem
3.6]{as4}, and if some $\gamma_i \neq 0$, then we are in the
situation of Proposition \ref{a2nequal3}. We note that
$A(\Gamma,V,\mu_1, \mu_2,\lambda,\gamma_1,\gamma_2)=
M_{\mu_1}/<x^n_2 - \mu_2 ( g_2^n-1), \upsilon -
\lambda(g_1^ng_2^n -1)>$.

Since $M_{\mu_1}$ is obtained by adjoining $x_1, z$ and $x_2$ via
Ore extensions to $k\Gamma$ and then factoring by a Hopf ideal,
then we may let $K$ be the Hopf subalgebra of $M_{\mu_1}$
generated by $\Gamma'$, the subgroup of $\Gamma$ generated by
$g_1$ and $g_2$, and by $x_2^n$ and $z^n$, i.e. by $g_1, g_2,
x_2^n $ and $\upsilon$.  Since $b_{ij}^n =1$, the $g_i$ commute
with $x_2^n$ and $ z^n$. Also $z^n$ and $x_2^n$ commute. For, if
$\gamma_2 =0$, then $x_2z = b_{21} zx_2$ and since $b_{21}^n =1$,
the commutation is clear. If $\gamma_2 \neq 0$, then $n=3$, and
the commutation of $x_2^3$ and $z^3$ follows from Lemma
\ref{aslemma} with $Y =z, X= x_2, \alpha =q^2, Z= \gamma_2 (g_1
g_2^2-1), \beta=1$. Thus the Hopf algebra $K$ is a commutative
polynomial algebra over $k\Gamma'$ in the indeterminates $x_2^n$
and $z^n$.

Now let  $\psi : K \to k$ be the algebra map defined by $\psi
(g_1)= \psi(g_2) =1$  , $\psi (x_2^n) = \mu_2$ and
$\psi(\upsilon) =  \lambda$. Then $\psi^{-1} (x_2^n) = -\mu_2$
and $\psi^{-1} (\upsilon) = -\lambda$. By Lemma \ref{conjugate},
the ideal $J$ generated by the skew-primitives $x^n_2$ and
$\upsilon$ and the ideal $I$ generated by the skew-primitives
$x^n_2 - \mu_2 (g_2^n-1 )$ and $\upsilon - \lambda (g_1^n g_2^n
-1)$ are conjugate in $K$. Also $(\psi J) \neq M_{\mu_1}$ since
$(\psi J)$ is the Hopf ideal generated by $x_2^n + \mu_2g_2^n$
and $v + \lambda g_1^ng_2^n$. Thus $M_{\mu_1}/(J) =
A(\Gamma,V,\mu_1, 0,0,\gamma_1,\gamma_2)$ and $M_{\mu_1}/(I) =
A(\Gamma,V,\mu_1, \mu_2,\lambda,\gamma_1,\gamma_2)$ are
quasi-isomorphic.

  Next, let $M = {\cal U} / < x^n_2,\; zx_1- b_{21} x_1z -
 \gamma_1 (g_1^2 g_2 -1), \; zx_2 - b_{12} b_{22} x_2 z -
 \gamma_2 (g_1g^2_2 -1)>$. Then $M/<x^n_1, \upsilon> \cong A(\Gamma,V,0,0,0,
 \gamma_1, \gamma_2)$ and $M/ <x_1^n -(g_1^n -1), \upsilon> \cong A(\Gamma,V,
 1, 0,0,
 \gamma_1, \gamma_2)$ and showing that
$J =< x_1^n, \upsilon>$ and
 $I =< x_1^n -( g^n_1 -1), \upsilon>$ are conjugate  Hopf ideals in
some Hopf subalgebra of $M$ will
 complete the proof. Let $K$ be the Hopf subalgebra of $M$ generated
 by $g_1, g_2, x_1^n$ and $z^n$. Again, $g_1$ and $g_2$ commute with
 $x_1^n$ and $z^n$ and if $\gamma_1 =0$, then $zx_1 = b_{21} x_1z$
 in $M$ so that $x_1^n$ and $z^n$ commute. If $\gamma_1 \neq 0$,
 then the commutation of $x_1^3$ and $z^3$ again follows from
 Lemma \ref{aslemma}. Define an algebra map $\varphi:K \to k$
 by $\varphi (g_1) = \varphi(g_2) = 1$ for $h \in \Gamma$, $\varphi (x_1^n) = 1$,
 $\varphi (v) =0.$ Then Lemma \ref{conjugate} again yields that
 $J$ and $I$ are conjugate in $K $ and so
 $A(\Gamma,V, \mu_1, \mu_2,\lambda, \gamma_1, \gamma_2) \sim
  A(\Gamma,V,0,0,0, \gamma_1, \gamma_2).$  If $n>3$ or $\Gamma$ is
  cyclic of order 3, then $\gamma_1=\gamma_2=0$, and all liftings
  of ${\cal B}(V) \# k\Gamma$ are quasi-isomorphic.

  Assume now that $n=3$ and that we are in the situation of
Proposition \ref{a2nequal3} with $\chi_1^2\chi_2 =
\chi_1\chi_2^2=\epsilon$.  Let $L= {\cal
U}/<x_1^3,x_2^3,zx_1-qx_1z,z^3>$, and let $K$ be the commutative
subHopf algebra of $L$ generated by $g_1,g_2$ and the
skew-primitive $zx_2-q^2x_2z$. Define an algebra map $\varphi:K
\to k$ by $\varphi(g_1)=\varphi(g_2)=1$ and
$\varphi(zx_2-q^2x_2z)= \gamma_2$. Then as above , the Hopf
ideals $J$ generated by $zx_2-q^2x_2z$ and $I$ generated by
$zx_2-q^2x_2z -\gamma_2 (g_1g_2^2 -1)$
 are conjugate in $S$ and so $L/(J)$ and $L/(I)$ are
 quasi-isomorphic, i.e. $A(\Gamma,V, 0,0,0,0,0) \sim
 A(\Gamma,V, 0,0,0,0,\gamma_2)$.  Similarly, $A(\Gamma,V, 0,0,0,0,0) \sim
 A(\Gamma,V, 0,0,0,\gamma_1,0).$ \qed

{\bf Question:}  For $n=3$ and $\gamma_1, \gamma_2$ nonzero, is
$A(\Gamma,V,0,0,0,\gamma_1,\gamma_2)\sim A(\Gamma,V,0,0,0,0,0)$?

{\bf Added in proof:} A. Masuoka has answered this question in the
affirmative. His method of proof is very much in the style of the
proofs in \cite{mas}.

Now we consider the case where $V$ is of type $B_2$ and $n \neq
5$.  If $A \cong U/J$ is the lifting determined by the scalars
$\mu_1, \mu_2, \lambda, \gamma$ as in Theorem \ref{quotient},
then we write $A=A(\Gamma,V,\mu_1,\mu_2,\lambda,\gamma).$

\begin{te}\label{quasib2}
For $V \in {^{k\Gamma}_{k\Gamma}{\cal YD}}$ of type $B_2$ and $n \neq
5$, any two liftings of ${\cal B}(V) \# k\Gamma$ are
quasi-isomorphic.
\end{te}

{\bf Proof.}  As in the proof of Theorem \ref{quasia2}, if
$b_{ij}^n \neq 1$, then only the lifting $A(\Gamma,V,0,0,0,0)$ is possible.
Therefore we assume that $b_{12}^n =b_{21}^n=1$.

We first show that
$A(\Gamma,V,0,\mu_2,0,0) \sim A(\Gamma,V,\mu_1, \mu_2, \lambda, \gamma)$
for fixed $\mu_2$ and any $ \mu_1, \lambda, \gamma$. Let
 $M(\mu_2) =U /< x_2^n -\mu_2 (g_2^n -1)>$, where $U$ is defined just before
Theorem \ref{quotient}.

Recall that $\upsilon, \omega$ defined in Section 2, equations
(\ref{upsilon}), (\ref{omega}), are skew-primitives, as is
$x_1^n$. We write $x_1^n, \upsilon, \omega$ also for the images
of these elements in $M(\mu_2)$ and note that they are still
skew-primitive.

Now let $K$ be the Hopf subalgebra of $M(\mu_2)$ generated by
  by $g_1$, $g_2$,
$x_1^n$, $z^n$ and $u^n$. As in the proof of Theorem
\ref{quasia2}, as an algebra $K$ is a commutative polynomial
algebra over $K\Gamma '$ where $\Gamma '$ is generated by $g_1,
g_2$.

Let $\psi$ be the $k$-linear map from $K$ to $k$ defined by
$\psi(g_1) = \psi(g_2) = 1  , \quad \psi(x^n_1) = \mu_1, \quad
\psi (\upsilon) = \lambda,\quad \psi (\omega) = \gamma.$ Then
$\psi$ defines an algebra map from $K$ to $k$. In $K$, let $J$ be
the Hopf ideal generated by $x_1^n, \upsilon$ and $\omega$ and
let $L$ be the Hopf ideal generated by the skew-primitives $x_1^n
- \mu_1 (g_1^n-1), \upsilon - \lambda (g_1^n g_2^n-1)\mbox{ and }
\omega - \gamma(g_1^n g_2^{2n} -1).$

Then by Lemma \ref{conjugate}, $J$ and $L$ are conjugate in $K$
and by \cite{mas}, $A(\Gamma,V,0,\mu_2, 0,0) \cong M(\mu_2)/(J)
\sim M(\mu_2)/(L) \cong A(\Gamma,V,\mu_1,\mu_2, \lambda,\gamma).$

Finally we show that $A(\Gamma,V,0,0,0,0) \sim
A(\Gamma,V,0,1,0,0).$ Let $M=U/<x_1^n>$, let $K$ be the
commutative Hopf subalgebra of $M$ generated by
$g_1,g_2,x_2^n,z^n$ and $u^n$, let $\psi:K \rightarrow k$ be the
algebra map defined by $\psi(g_i)=1= \psi(x_2^n), \psi(\upsilon)
=\psi(\omega)=0.$ Now the same argument finishes the proof.  \qed

\appendix
\section{A generalization of the $q$-binomial theorem}
{\bf  \author{by M. Beattie, S. D\u{a}sc\u{a}lescu, \c{S}. Raianu
and I. Rutherford \footnote{Ian Rutherford was supported by an
NSERC Undergraduate Student Research Award at Mount Allison U. in
the summer of 2000.}}}

Throughout, we work over a field $k$, not necessarily
algebraically closed.

>From Theorem \ref{classic} (i), it is straightforward to prove
that:
\begin{eqnarray}\label{qpasc}
{n+1 \choose i}_q {i \choose j}_q={n \choose i-1}_q {i-1 \choose j-1}_q
+ q^i{n \choose i}_q {i \choose j}_q + q^{j}{n \choose i-1}_q
{i-1 \choose j}_q.
\end{eqnarray}

Now we prove the generalized quantum binomial theorem used in the
calculations in this paper.  However, this theorem is interesting
in its own right.

\begin{te}\label{main}  Suppose that $q \in k^*$ and $\lambda \in k$
and for $x,t,z $
in some $k$-algebra, we have the following relations :
\begin{eqnarray}\label{assumptions}
x z = q z x;\mbox{   } z t = q t z; \mbox{  }x t = qt x + \lambda z^2.
\end{eqnarray}
Then
\begin{displaymath}
(x + bz + t)^n =
\sum_{i=0}^n\sum_{j=0}^i{ n \choose i}_{q}{i \choose j}_q
\nu(i-j)t^jz^{i-j}x^{n-i}
\end{displaymath}
where $\nu = \nu_{b,\lambda}$ is a function from {\bf N} to $k$
defined recursively
by $\nu(0) = 1, \nu(1) = b$ ,  and
$ \nu(n) = b\nu(n-1) + \lambda (n-1)_{q}\nu(n-2), \mbox{ for }n \geq 2.$\\

\end{te}

{\bf Proof.}
The proof is by induction. The formula can easily be checked
for $n = 1,2$.
Now assume
that the formula holds for $n=k$, and we show that it is valid for $n=k+1$.
First we note that it follows directly from Lemma \ref{aslemma}
or from a simple induction
argument that
\begin{eqnarray} \label{xt}
xt^n = q^{n}t^nx + \lambda q^{n-1}(n)_{q}t^{n-1}z^2.
\end{eqnarray}

Then we compute
\begin{eqnarray*}
(x + bz + t )^{k+1} &=& (x+bz+t)(x+bz+t)^k
\\
&=&  \sum_{i=0}^k \sum_{j=0}^i {k \choose i}_q {i \choose j}_q
\nu(i-j)(x+bz+t)t^jz^{i-j}x^{k-i}
\\
&=& \sum_{i=0}^k\sum_{j=0}^i {k \choose i}_q {i \choose j}_q
\nu(i-j)((q^{j}t^jx
+ \lambda q^{j-1}(j)_{q}t^{j-1}z^2 )z^{i-j}x^{k-i}
\\
&+&  q^{j}bt^jz^{i-j+1}x^{k-i}
+ t^{j+1}z^{i-j}x^{k-i}) \mbox{ by ( \ref{xt})}
\\
&=& \sum_{i=0}^k\sum_{j=0}^i {k \choose i}_q {i \choose j}_q
\nu(i-j)(\lambda q^{j-1}(j)_{q}t^{j-1}z^{2+i-j}x^{k-i}
\\
&+&  q^{i}t^jz^{i-j}x^{k+1-i}
+ q^{j}bt^jz^{i-j+1}x^{k-i}
+ t^{j+1}z^{i-j}x^{k-i}).
\end{eqnarray*}

Now let $i'=i+1,j'=j-1, j''=j+1$ and then $(x + bz + t)^{k+1}=$
\begin{eqnarray*}
&=& \sum_{i'=1}^{k+1}\sum_{j'=-1}^{i'-2}{k \choose i' -1}_q
{{i'-1} \choose {j'+1}}_q
\nu(i'-j'-2)\lambda q^{j'}(j'+1)_{q}t^{j'}z^{i'-j'}x^{k+1-i'}
\\
&+& \sum_{i=0}^k\sum_{j=0}^i {k \choose i}_q {i \choose j}_q
\nu(i-j)q^{i}t^jz^{i-j}x^{k+1-i}
\\
&+& \sum_{i'=1}^{k+1}\sum_{j=0}^{i'-1}{k \choose i'-1}_q {i'-1 \choose j}_q
\nu(i'-j-1)q^{j}bt^jz^{i'-j}x^{k+1-i'}
\\
&+& \sum_{i'=1}^{k+1}\sum_{j''=1}^{i'}{k \choose i' - 1}_q{i' -1
\choose j'' -1}_q
\nu(i'-j'')t^{j''}z^{i'-j''}x^{k+1-i'},
\end{eqnarray*}
and then, using the fact that ${m \choose s}_q = 0 \mbox{ if } s>m
\mbox{ or } s<0$ and
$(0)_q = 0$, we have that $(x + bz + t)^{k+1}=$

\begin{eqnarray*}
&=& \sum_{i=0}^{k+1}\sum_{j=0}^{i}{k \choose i -1}_q
{{i-1} \choose {j}}_q \frac{(i-j-1)_{q}}{(j+1)_{q}}
\nu(i-j-2)\lambda q^{j}(j+1)_{q}t^{j}z^{i-j}x^{k+1-i}
\\
&+& \sum_{i=0}^{k+1}\sum_{j=0}^i {k \choose i}_q {i \choose j}_q
\nu(i-j)q^{i}t^jz^{i-j}x^{k+1-i}
\\
&+& \sum_{i=0}^{k+1}\sum_{j=0}^{i}{k \choose i-1}_q {i-1 \choose j}_q
\nu(i-j-1)q^{j}bt^jz^{i-j}x^{k+1-i}
\\
&+& \sum_{i=0}^{k+1}\sum_{j=0}^{i}{k \choose i - 1}_q{i -1 \choose j -1}_q
\nu(i-j)t^{j}z^{i-j}x^{k+1-i}
\\
&=& \sum_{i=0}^{k+1}\sum_{j=0}^{i}({k \choose i - 1}_q{i -1 \choose j}_q
(b \nu(i-j-1) + \lambda \nu(i-j-2)(i-j-1)_q)q^j
\\
&+& {k \choose i }_q{i  \choose j}_q \nu(i-j)q^i
+{k \choose i - 1}_q{i -1 \choose j-1}_q \nu(i-j))t^jz^{i-j}x^{k+1-i}
\\
&=& \sum_{i=0}^{k+1}\sum_{j=0}^{i}({k \choose i - 1}_q{i -1 \choose j-1}_q
 + q^i {k \choose i}_q{i \choose j}_q +
q^j {k \choose i - 1}_q{i -1 \choose j}_q ) \nu(i-j) t^jz^{i-j}x^{k+1-i}
\\
&& \mbox{      by the definition of } \nu
\\&=& \sum_{i=0}^{k+1}\sum_{j=0}^i{ k+1 \choose i}_{q}{i \choose j}_q
\nu(i-j)t^jz^{i-j}x^{k+1-i} \mbox{ by (\ref{qpasc}), as required. \qed }
\end{eqnarray*}

 \begin{res} \label{mainremarks}
{\rm {\bf i.} If $\lambda = 0 $, then $(x + bz + t)^n = ((x+bz)+t)^n$
where $(x+bz)t=qt(x+bz) $ and so the same result may be obtained
directly from the $q$-binomial theorem (Theorem \ref{classic}(ii)).
\\
{\bf ii.} Suppose $q$ is a primitive $n$th root of unity.
Then ${n \choose i}_q=0$
unless $i=0$ or $i=n$ and the formula in Theorem \ref{main}
becomes $(x + bz + t)^n = x^n + \nu_{b,\lambda}(n)z^n + t^n.$
}
\end{res}

The description of $(x + bz + t)^n $ would now be complete if we
had a general
formula for $\nu_{b,\lambda}(s).$

\begin{pr}
If $q=1$ and $b \ne 0$ then for $n \ge 0$,
$$ \nu(n) = \sum _{i=0}^{n/2} {n \choose 2i} \frac{(2i)!}{2^i (i!)}
b^{n-2i} \lambda^i$$
\end{pr}
{\bf Proof.} Note that since ${ n \choose 2i}=0$ if $2i>n$, the summation
is from $0$ to $\lfloor n/2 \rfloor$ .  First, if $n=0, \mbox{ we check that }
{0 \choose 0}
\frac{0!}{2^0 0!}b^0 \lambda ^0 = 1$ and , if $n=1$, the formula gives
$ { 1 \choose 0} \frac {0!}{2^0 0!}b^1 = b$ .  Now suppose the formula holds
for $n \le k+1$ and we compute $\nu(k+2) = b \nu(k+1) + (k+1)
\lambda \nu(k) $.
By the induction assumption, this is
\begin{eqnarray*}
&& \sum _{i=0}^{(k+1)/2} {k+1 \choose 2i} \frac{(2i)!}{2^i (i!)}
b^{k+2-2i} \lambda^i
 + \sum _{i=0}^{k/2} {k \choose 2i}(k+1) \frac{(2i)!}{2^i (i!)} b^{k-2i}
 \lambda^{i+1}\\
&=& \sum _{i=0}^{(k+2)/2} {k+1 \choose 2i} \frac{(2i)!}{2^i (i!)}
b^{k+2-2i} \lambda^i
 + \sum _{i'=0}^{(k+2)/2} {k \choose 2i'-2}(k+1)
\frac{(2i'-2)!}{2^{i'-1} ((i'-1)!)} b^{k+2-2i'} \lambda^{i'}  \\
&& \mbox { where } i'=i+1 \\
&=& \sum _{i=0}^{(k+2)/2}\big( {k+1 \choose 2i}\frac{(2i)!}{2^i (i!)} +
{k+1 \choose 2i-1}\frac{(2i-1)!}{2^{i-1} ((i-1)!)}\big) b^{k+2-2i}\lambda^i\\
&=& \sum _{i=0}^{(k+2)/2}{ k+2 \choose 2i} \frac {(2i)!}{2^i i!}
b^{n+2-2i}\lambda^i
\mbox{ .             \qed}
\end{eqnarray*}

If $q \neq 1,$ there is a formula for the computation of $\nu(n)$
in terms of $\alpha, \beta$ ($\alpha, \beta$ possibly lie in some
extension field of $k$) where $\alpha + \beta =b$ and $\alpha
\beta = \lambda/(q-1).$

\begin{pr}\label{nu} Suppose that $q \neq 1.$
Let $\alpha, \beta$ be the roots of the polynomial $Y^2 -bY + \lambda/(q-1)$
in some extension field $ \overline{k}$ of $k$. Then
$$ \nu_{b,\lambda}(n) = \sum_{i=0}^{n} {n \choose i}_q \beta^i
\alpha^{n-i}.$$
\end{pr}
{\bf Proof.} Since $x + bz + t=(x + \alpha z) + (\beta z + t)$,
and since $(x + \alpha z)(\beta z + t) = q(\beta z + t)(x +
\alpha z)$, we have by the $q$-binomial theorem (Theorem
\ref{classic}) that
\begin{eqnarray*}
(x + bz + t)^n &=& \sum_{i=0}^n{n \choose i}_q(\beta z + t)^i(x +
\alpha z)^{n-i}
\\
&=& \sum_{i=0}^n{n \choose i}_q (\sum_{j=0}^i{i \choose j}_qt^j
(\beta z)^{i-j})
(\sum_{k=0}^{n-i}{n-i \choose k}_q(\alpha z)^kx^{n-i-k}).
\end{eqnarray*}

Comparing terms in this formula with the one in Theorem \ref{main},
we see that
\begin{eqnarray*}
\sum_{m=0}^n {n \choose m}_q \nu(m) z^m x^{n-m} &=&
\sum_{i=0}^n \sum_{k=0}^{n-i}{n \choose i}_q
{n-i  \choose k}_q(\beta z)^i(\alpha z)^kx^{n-i-k}
\\
&=& \sum_{w=0}^n \sum_{i=0}^w {n \choose i}_q {n-i \choose w-i}_q
\beta^i \alpha ^{w-i}z^wx^{n-w}.
\end{eqnarray*}

Comparing coefficients of $z^m x^{n-m}$, we obtain
$${n \choose m}_q\nu(m) = \sum_{i=0}^m {n \choose i}_q{n-i \choose m-i}_q
\beta^i \alpha^{m-i},$$
and then letting $m=n$, we obtain the statement.             \qed

\begin{co}\label{conu}
If $q$ is a primitive $n$th root of unity, and $n>1$,
then $\nu(n)= \alpha^n + \beta^n \in k$.
\end{co}

If $b$=0, then $\nu(1)=b=0$, and since $\nu(2n+1)=
b\nu(2n) + \lambda(2n)_q \nu (2n - 1)$, it is clear that $\nu(2n+1)=0$
for $n \geq 0$. However, $\nu(0)=1, \nu(2) = \lambda(1)_q ,
\nu(4)= \lambda^2 (3)_q$, and , in general, $\nu(2n)=
\lambda^n (2n-1)_q (2n-3)_q \ldots (1)_q$.

\begin{co}(to Theorem \ref{main})
For $x,z,t,q,\lambda$ satisfying (\ref{assumptions}),
$$(x+t)^n = \sum_{i=0}^n \sum_{m=0}^{\lfloor\frac{i}{2}\rfloor}
{n \choose i}_q
{i \choose i-2m}_q \nu(2m)t^{i-2m}z^{2m}x^{n-i}.$$
\end{co}
{\bf Proof.} Applying Theorem \ref{main}
with $b=0$, we have
\begin{eqnarray*}
(x + t)^n &=&
\sum_{i=0}^n\sum_{j=0}^i{ n \choose i}_{q}{i \choose j}_q
\nu(i-j)t^jz^{i-j}x^{n-i}
\\
&=& \sum_{i=0}^n\sum_{k=0}^i{ n \choose i}_{q}{i \choose i-k}_q
\nu(k)t^{i-k}z^{k}x^{n-i}
\\
&=& \sum_{i=0}^n\sum_{m=0}^{\lfloor\frac{i}{2}\rfloor}
{ n \choose i}_{q}{i \choose i-2m}_q
\nu(2m)t^{i-2m}z^{2m}x^{n-i}.  \mbox{\qed}
\end{eqnarray*}

\begin{re}\label{serban}
{\rm Suppose $x,t,s$ are such that
$$xs = q^2sx, \; st =q^2ts \; \mbox{ and } xt =qtx + \lambda s,$$
i.e. the relations (\ref{assumptions}) hold with $s=z^2$.  Then if $q$
is a primitive $n$th root of unity and $n$ is odd, $(x+t)^n = x^n + t^n.$}
\end{re}

\begin{center}
{\bf Acknowledgment}\\
\end{center}
Thanks to N. Andruskiewitsch for several interesting
conversations about the lifting problem.  Also many thanks to A.
Masuoka for his interest in this paper and for his many helpful
comments on the problem of whether liftings are quasi-isomorphic,
including providing an affirmative answer to the question of
whether all liftings of type $A_2$ are quasi-isomorphic.

{\small

\begin{flushleft}

\author{M. Beattie and I. Rutherford, \\Department of Mathematics\\
and Computer Science\\
Mount Allison University\\Sackville, N.B., Canada E4L 1E6\\
\vspace{12pt}
\and S. D\u{a}sc\u{a}lescu, \\
University of Bucharest\\
Faculty of Mathematics, \\Str. Academiei 14\\ RO-70109 Bucharest 1,
Romania  \\
\vspace{12pt}
\and \c{S}. Raianu\\
Department of Mathematics\\
Syracuse University\\
Syracuse, New York 13244, U.S.A.\\}

\end{flushleft}
 }

\end{document}